\documentclass{article}
\usepackage{graphicx}
\usepackage{amssymb,amsfonts,amsthm}
\usepackage{amsmath,amsthm,amssymb}
\usepackage{amsfonts}

\newtheorem{theorem}{Theorem}[section]
\newtheorem{lemma}[theorem]{Lemma}

\theoremstyle{definition}

\newtheorem{rk}[theorem]{Remark}

\newcounter{ppp}

\newcommand{\la}{\langle}
\newcommand{\ra}{\rangle}

\begin{document}

\renewcommand{\theequation}{\thesection.\arabic{equation}}

\title{Linear automorphism groups of relatively free groups}
 \author{ A.Yu. Olshanskii\thanks{\large The
author was supported in part by the NSF grants DMS
0455881, DMS 0700811, and by the Russian Fund for Basic Research grant
05-01-00895.}}
\date{}
\maketitle



\section{Introduction}

Let $F=F(x_1,x_2,\dots)$ be an absolutely free group with basis $x_1,x_2,\dots$. 
Recall that a group $H$ satisfies an identity $w(x_1,\dots,x_n)=1$ for a word $w(x_1,\dots,x_n)=w\in F$ if
 $w$ vanishes under every homomorphism $F\to H$. A {\it variety} of groups
 is a class of groups consisting of all groups that satisfy some set of identities. 
 
 A {\it relatively free group} $G$ is a free group in a group variety $\cal V$, i.e.,
 $G$ belongs to $\cal V$, and $G$ is generated by a set $Y$ such that every mapping
 $Y\to H$, where $H\in \cal V$, extends to a homomorphism $G\to H$. The (free) rank
 of $G$ is the cardinality of $Y$. (See details in \cite{N}.) We denote by $F_n$,
 by $A_n=F_n/F'_n$, and by $M_n=F_n/F''_n$ the absolutely free group, the free abelian and the free metabelian groups of rank $n$,
 respectively. 
 
 The group of inner automorphisms of a group $G$ is normal in $Aut(G)$, and so
 the factor group $G/Z(G)$, where $Z(G)$ is the center of $G$, can be canonically
 identified with a normal subgroup of $Aut(G)$.
 Since the center of $M_n$ is trivial for $n\ge 2$ (see \cite{N}, 25.63), one can identify $M_n$
 with the normal subgroup of $Aut(M_n)$ consisting of the inner automorphisms of $M_n$.
 
 Our paper is inspired by the following result of V.P.Platonov that answers a question raised 
 by H.Mochizuki \cite{M}.
 
 \begin{theorem} (V.P.Platonov \cite{P}) \label{plat} (1) Let $\rho$ be a finite-dimensional 
 linear representation of the automorphism group $Aut(M_n)$ over a field $k$. Then the 
 image $\rho(M_n)$ is a virtually nilpotent group. (2) It follows that the group $Aut(M_n)$
 is not linear for $n>1$. 
 \end{theorem}
 
 As usual, a group $G$ is called {\it linear} if it is isomorphic to a subgroup of $GL_m(k)$
 for some field $k$ and some integer $m\ge 1$.
 
 In Section 2, we give an
 alternative and shorter proof of Theorem \ref{plat}.  
 Then a similar approach and the utilization of some known properties of group varieties lead to the complete 
 description of relatively free groups $G$ for which $Aut(G)$ is a linear group. 
 
 \begin{theorem}\label{2} Let $G$ be a relatively free  but not absolutely free group. 
 The automorphism group $Aut(G)$ is linear if and only if $G$ is a finitely generated 
 virtually nilpotent group. 
 
 Furthermore, if the group $G$ is  finitely generated but not virtually-nilpotent, 
 then there is an automorphism $\phi$ of $G$ such
 that the extension $P$ of $G/Z(G)$ by $\phi$ is a non-linear subgroup of $Aut(G)$; and if 
 $G$ is finitely generated and virtually nilpotent, then the holomorph $Hol(G)$ is linear over $\mathbb Z$.
\end{theorem}

Recall that a group $G$ is {\it virtually nilpotent} if it contains a (normal) nilpotent
subgroup of finite index. (The "if" part of the statement does not need the
hypothesis that $G$ is relatively free.)

\begin{rk} The automorphism group $Aut(F_n)$ is not linear for $n\ge 3$ (Formanek, Procesi \cite{FP})
but the group $Aut(F_2)$ is linear (Krammer \cite{Kr}).

\end{rk} 

\begin{rk} Note that the formulation of Theorem \ref{2} is similar to those contained in
the papers of O.M.Mateiko and O.I.Tavgen' \cite{MT} and A.A.Korobov \cite{Ko}. Nevertheless  we prove Theorem \ref{2} here
for the following reasons. 
(1) There is a mistake in both \cite{MT} and \cite{Ko}. Namely, the proofs essentially use
the "known property" of Fitting subgroups to be fully characteristic. But this
does not hold even for relatively free groups. (For example, the Fitting subgroup is not fully characteristic 
in the free group of rank $n>1$ of the variety generated by the alternating group $Alt(5)$.) (2) The formulation
of the main theorems is not quite correct in both \cite{MT} and \cite{Ko} because it is not proved there that the 
virtually nilpotency of a free group of rank
$n>1$ in a variety implies virtually nilpotency of free groups having rank $>n$. (3) Our
proof is simpler. (The authors of both papers \cite{MT} and \cite{Ko} refer, in particular, to the statement 
that all locally finite groups of exponent dividing $m$ form a variety. This claim is equivalent to the restricted
Burnside problem for groups of exponent $m$, and the affirmative solution is bases on the Classification Hypothesis for
finite simple groups.)
\end{rk}
 
 \section{Free matabelian case}
 
 The following Kolchin--Mal'cev Theorem is the most--known fact on
 linear solvable groups.
 
 \begin{lemma}\label{linsol} Every linear solvable group has a subgroup $H$ of
 finite index such that the derived subgroup $H'$ is nilpotent. $\Box$ 
 \end{lemma}\endproof
 
 Let $\phi$ be an automorphism of the free abelian group $A_n$ and $B_{\phi}=\la A_n,\phi\ra$
 the extension of $A_n$ by the automorphism $\phi$. Assume that no root of the characteristic 
 polynomial of $\phi$ is an $k$-th root of $1$ for any integer $k>0$. The following property
 of the group $B_{\phi}$ is folklore.
 
 \begin{lemma}\label{phi} Let $C$ be a subgroup of finite index in $B_{\phi}$. Then the derived
 subgroup $C'$ has finite index in $A_n$.
 \end{lemma}
 \proof
 Since $C$ is of finite index in $B_{\phi}$, it must contain $\phi^k$ and $mA_n$ for some 
 positive integers $m$ and $k$. (We use the additive notation for $A_n$ in this proof.)
 Therefore $C'$ contains the subgroup $[m A_n,\phi^k]=\{\phi^k(ma)-ma\mid a\in A_n\}.$ 
 
  Proving by contradiction, assume that the index $(A_n:[m A_n,\phi^k])$ is infinite. Then the 
 image of $mA_n$ under the mapping $\phi^k-id$ is of rank $<n$. Hence $0$ is an eigenvalue of $\phi^k-id$,
 and so $\lambda^k=1$ for a characteristic root $\lambda$ of $\phi$; a contradiction. \endproof
 
 {\it Proof of Theorem \ref{plat}.} 
 There is nothing to prove if $n=1$. For every $n\ge 2$, there exists an automorphism $\phi$ of $A_n$ whose characteristic
 roots are not roots of $1$. (One can easily find such automorphisms for $n=2,3$ and note for $n>3$ that $A_n$ 
 is a direct sum of subgroups isomorphic to $A_2$ and $A_3$.) We keep the same notation  $\phi$ for a lifting of $\phi$ to $Aut(M_n)$.
 (Recall that $A_n\simeq M_n/M'_n$, and the induced homomorphism $Aut(M_n)\to Aut(A_n)$ is surjective by \cite{N}, 41.21.) 
   
 Let $P =\la M_n,\phi\ra$ be the extension of $M_n$ by the automorphism $\phi$. This $P$ is a 
 subgroup of $Aut(M_n)$ since $\la\phi\ra\cap M_n=\{1\}$; and $P$ is solvable 
 because the factor-group $P/M_n$ is cyclic. By Lemma \ref{linsol}, there is a normal
 subgroup $T$ of finite index in $P$ such that the subgroup $\rho(T')$ is nilpotent. 
 
 The canonical image $C=TM'_n/M'_n$ of $T$ in $B_{\phi}=P/M'_n$ has finite index, 
 and therefore $C'$ is of finite index in $A_n=M_n/M'_n$ by Lemma \ref{phi}. Hence the inverse image $D=T'M'_n$ 
 is of finite index in $M_n$. 
 
 The normal subgroup
 $\rho(M'_n)$ of $\rho(M_n)$ is abelian since $M_n$ is metabelian. Thus, $\rho(D)$ is
 nilpotent being a product of two nilpotent normal in $\rho(P)$ subgroups $\rho(T')$ and $\rho(M'_n)$.
 Since $(M_n:D)<\infty$, the theorem is proved.
 $\Box$
 
 \section{Few lemmas on varieties of groups}
 
The {\it product}
${\cal U}{\cal V}$ of two group varieties contains all the groups $G$ having a normal
subgroup $N$ such that $N\in \cal U$ and $G/N\in \cal V$; ${\cal UV}$ is also a group
variety (\cite{N}, 21.12).

\begin{lemma} (\cite {Sh}).\label{shm} Let $L$ be a free group of rank $n\ge 2$ in a product
of varieties ${\cal UV}$. Assume that the free group of rank $n$ in the variety $\cal V$ is
infinite and ${\cal U}$ contains a non-trivial group. 
Then the center of $L$ is trivial. $\Box$ \end{lemma}

We denote by $\cal A$ (by ${\cal A}_k$) the variety
of all abelian groups (of all abelian groups of exponent dividing $k$),  
and denote by $M_{k,n}$ the free group of rank $n$ in
the variety  ${\cal M}_k={\cal A}_k\cal A$.

\begin{lemma} \label{AkA} If $n,k\ge 2$, then the group $M_{k,n}$ is not virtually nilpotent. 
\end{lemma} 

\proof The wreath product $W={\mathbb Z}_k wr {\mathbb Z}$ of a cyclic group of order $k$ and an infinite
cyclic group is 2-generated, and it belongs to the variety ${\cal M}_k$. Therefore $W$ is a homomorphic
image of the group $M_{k,n}$. Since $W$ is not virtually nilpotent, $M_{k,n}$ is not virtually nilpotent too. 
\endproof

A variety $\cal V$ is called {\it solvable} if all the groups of $\cal V$ are solvable.
 
\begin{lemma} (\cite{G}). \label{altern} Let $\cal S$ be a solvable variety of groups.
Then either $\cal S$ contains as a subvariety the product ${\cal M}_p={\cal A}_p\cal A$ for some prime $p$ 
or every finitely generated group in $\cal S$ is virtually nilpotent. $\Box$
\end{lemma}

Further we call a variety $\cal V$ {\it proper} if it does not contain all groups, or equivently,
the absolutely free group $F_2$ does not belong to $\cal V$. It is easy to see that a
product of two proper varieties is proper. (See also \cite{N}.)
The minimal variety containing a group $Q$ is denoted by $var Q$.
 Given a group $G$ and a variety $\cal V$, the {\it verbal subgroup} $V(G)$ corresponding to $\cal V$ is
 the smallest normal subgroup $N$ of $G$ such that $G/N\in \cal V$.

\section{Proof of Theorem \ref{2}}

By Auslander -- Baumslag's theorem \cite{AB}, the holomorph of any finitely generated (virtually) nilpotent group $G$ is linear over $\mathbb Z$. Thus, it remains to consider a non-virtually-nilpotent free group $G$ of rank $n\ge 2$ in a proper variety $\cal V$ and
construct the required automorphism $\phi$. (A non-trivial relatively
free group of infinite rank admits the automorphisms from an infinite symmetric group that is not linear.)
Now the quotient $H=G/Z(G)$ is non-virtually-nilpotent normal subgroup of $Aut(G)$. We may assume that $H$ is a linear group since
otherwise there is nothing to prove.

Since both $G$ and $H$ satisfy a non-trivial
identity and $H$ is linear, the group $H$ is virtually solvable by Platonov's theorem \cite{P1}, 
i.e., the solvable radical $R$ of $H$ is of finite index in $H$. Therefore $R$ is a finitely generated
but not virtually nilpotent solvable group.

By Lemma \ref{altern}, there is a prime $p$ such that ${\cal M}_p\subseteq var R\subseteq var H\subseteq var G.$

Therefore there are canonical epimorphisms $G\to M_{p,n}$ and $G\to A_n$. 
The kernels are the verbal subgroups of $G$ corresponding to the varieties ${\cal M}_p$ and $\cal A$, respectively. 
The latest kernel is just the derived subgroup $G'$, and we denote by $M_p(G)$ the former one.

The center of $M_{p,n}$ is trivial by Lemma \ref{shm}, and so $Z(G)$ is contained in  $M_p(G)\subseteq G'$. 
Consequently, we have isomorphisms $G/M_p(G)\simeq H/M_p(H)\simeq M_{p,n}$ and $G/G'\simeq H/H'\simeq A_n$. 

Now, as in the proof of Theorem \ref{plat}, we introduce an automorphism $\phi$ of $A_n$ whose action 
has no characteristic roots equal to any root of $1$. As there (by \cite{N}, 41.21) one can lift $\phi$ to $Aut(G)$ and also
to $Aut(H)$ since the center $Z(G)$ is a characteristic subgroup of $G$. Denote by $P= \la H,\phi\ra$ the
extension of $H$ by the automorphism $\phi$. It is a subgroup of $Aut(G)$ as in the proof of Theorem \ref{plat}.  

Proving by contradiction, assume that $P$ is a linear group. Since $P\in {\cal V}\cal A$, the group $P$
satisfies a non-trivial identity, and by \cite{P1}, $P$ must have a solvable normal subgroup of finite index. 
By Lemma \ref{linsol}, $P$ contains a normal subgroup $T$ of finite index with nilpotent derived subgroup $T'$.
Applying Lemma \ref{phi} to the image of $T$ in $B_{\phi}=P/H'$, we have $(H:(T'H'))<\infty$.  

The quotient $H'/M_p(H)$ is
an abelian normal subgroup of $H/M_p(H)\simeq M_{p,n}$. Since $T'$ is nilpotent and normal in $P$, 
the image of the subgroup $T'H'$
under the canonical epimorphism $H\to M_{p,n}$ is nilpotent too. But this image
is of finite index in $M_{p,n}$ because $(H:(T'H'))<\infty$. This contradicts 
the statement of Lemma \ref{AkA}, and so the theorem is proved. $\Box$

\bigskip

{\bf Acknowledgments.}
The author is grateful to V.D.Mazurov, V.A.Romankov, N.S.Romanovskiy,  and A.E.Zalesski for useful discussion.

\begin{minipage}[t]{3 in}
\noindent Alexander Yu. Olshanskii\\ Department of Mathematics\\
Vanderbilt University \\ alexander.olshanskiy@vanderbilt.edu\\
 and\\ Department of
Higher Algebra\\ MEHMAT,
 Moscow State University\\

\end{minipage}

\end{document}